\documentclass[english]{amsart}

\usepackage{babel}
\usepackage{amstext}
\usepackage{amsmath}
\usepackage{amsfonts}
\usepackage{amssymb}
\usepackage{latexsym}
\usepackage{ifthen}
\usepackage[all,2cell,dvips]{xy} \UseAllTwocells \SilentMatrices
\usepackage{xy}
\xyoption{all}
\newtheorem{lemma1}{}[section]

\newenvironment{lemma}{\begin{lemma1}{\bf Lemma.}}{\end{lemma1}}
\newenvironment{example}{\begin{lemma1}{\bf Example.}\rm}{\end{lemma1}}

\newenvironment{theorem}{\begin{lemma1}{\bf Theorem.}}{\end{lemma1}}
\newenvironment{proposition}{\begin{lemma1}{\bf Proposition.}}{\end{lemma1}}
\newenvironment{corollary}{\begin{lemma1}{\bf Corollary.}}{\end{lemma1}}
\newenvironment{remark}{\begin{lemma1}{\bf Remark.}\rm}{\end{lemma1}}
\newenvironment{definition}{\begin{lemma1}{\bf Definition.}}{\end{lemma1}}

\newenvironment{conjecture}{\begin {lemma1}{\bf Conjecture.}}{\end{lemma1}}

\makeatletter
\ifnum\@ptsize=0  \fi \ifnum\@ptsize=2  \fi 

\newcommand\sF{{\mathcal F}}
\newcommand\sG{{\mathcal G}}

\newcommand\sI{{\mathcal I}}

\newcommand\sO{{\mathcal O}}

\newcommand\bZ{{\mathbb Z}}

\newcommand\bQ{{\mathbb Q}}

\DeclareMathOperator*{\Sym}{Sym}

\title {BIRATIONAL STABILITY OF THE COTANGENT BUNDLE} 
\author{Fr{\'e}d{\'e}ric Campana}
\address{Fr\'ed\'eric Campana}
\email{campana@iecn.u-nancy.fr}
\date{\today}

\begin{document}

\begin{center}
{\Large }
\end{center}

\maketitle

{\bf Abstract:} Roughly, we show, using some standard conjectures, that the cotangent bundle of a complex projective manifold $X$ is `birationally' semi-stable, unless $X$ is uniruled, in which case the unstability is controlled by its `rational quotient'. More precisely, we introduce a birational invariant $\kappa_{++}(X)\geq \kappa(X)$, the difference measuring the birational unstability of its cotangent bundle. The invariant $\kappa_{++}(X)$ is the maximum of the  `saturated' Kodaira dimensions of rank-one coherent subsheaves of all $\Omega^p_X,$ for any $p>0$. This is a measure of the birational positivity of these subsheaves, in contrast to their `numerical' positivity, by means of polarisation-slopes. Conjecturally, $\kappa_{++}(X))=\kappa(R(X))$, $R(X)$ being the `rational quotient' of $X$(see \S 1). For example, one should have: $\kappa_{++}(X)=-\infty$ if and only if $X$ is rationally connected, and $X$ being unstable in this sense if and only if uniruled, but not rationally connected. The study of the Kodaira dimensions of such sheaves was initiated By F. Bogomolov in \cite{Bo78}, where bounds and a partial geometric description of extremal cases were established.

We extend these notions and conjectures to the category of `smooth orbifolds'. These appear naturally in the geometric interpretation of the saturation process of the subsheaves of $Sym^N(\Omega_X^p)$ introduced in the definition of $\kappa_{++}(X)$. This category is, on the other hand, needed in an essential way for the birational classification. We then reduce the above conjecture in the `orbifold' setting to other standard conjectures of the LMMP, and to an extension of Miyaoka's generic semi-positivity for lc\footnote{Stands for `log-canonical', see \cite{KM}, for example, for this notion, as well as for klt pairs.} pairs with $c_1=0$. 
The notion of rational curve, uniruledness and rational connectedness will be introduced in the context of smooth orbifods as well. We show the uniruledness of some peculiar Fano orbifolds by specific elementary methods.

We also prove a stronger `numerical dimension' version: $\nu_{++}(X)=\kappa_{++}(X)=\kappa(R(X))$ in the orbifold context, conditionally in the same set of conjectures in \S 7 (see definitions there). The proof simplifies the earlier proof given in \cite{Ca 09} of the weaker result concerning $\kappa_{++}$.

We incidentally give in \S 8 a (seemingly) new very simple proof of the pseudo-effectivity of the relative canonical bundle of a fibre space when its generic fibres are not uniruled. This gives a weakened version of Viehweg's weak positivity results, which permits to deduces $C_{n,m}$ conjecture from Abundance conjecture, and is potentially susceptible of further developments.

As an application outside of the birational classification, we will mention the isotriviality conjecture for families of canonically polarised manifolds parametrised by a `special' quasi-projective manifold, which can also be reduced also to the very same set of conjectures, and thus becomes a problem in birational classification.

The present text complements results in \cite{Ca 09}, \cite{Ca07}, and \cite{Ca04}, where complete definitions can be found.  

\

\tableofcontents

\section{\bf  Introduction.}

Let $X$ be a complex projective connected $n$-fold. The main concern of birational algebraic geometry consists in deducing qualitative birational geometric properties of $X$ from positivity or negativity properties of the canonical bundle $K_X$. In particular, one would like to describe in these terms the birational invariants of manifolds $Y$ which are `rationally dominated' by $X$, i.e: such that there exists a dominant rational map $f:X\dasharrow Y$ (we then write: $Y\leq X)$, and so the following invariant:.

\begin{definition}\label{dk+} Let $\kappa_+(X):=max_{\{Y\leq X\}}\{\kappa(Y)\}$.

Thus: $n\geq \kappa_+(X)\geq \kappa(X)$, and $\kappa_+(X)=-\infty$ iff $\kappa(Y)=-\infty $ for any $Y\leq X$. 
\end{definition}

This invariant has a conjectural description, in terms of the `rational quotient' (or `MRC-fibration' ) $r_X:X\dasharrow R(X)$. Recall that this rational fibration has rationally connected (R.C, for short) fibres, and non-uniruled base $R(X)$ (or is a point iff $X$ is rationally connected). When $X$ is non-uniruled, $R(X)=X$, and $r_X$ is just the identity map.

\begin{conjecture}\label{ck+} For any $X$, $\kappa_+(X)=\kappa(R(X))$. In particular, $\kappa_+(X)=\kappa(X)$ if $\kappa(X)\geq 0$.
\end{conjecture}

This conjecture can be reduced to two other quite standard conjectures. Recall that a (rational) `fibration' means here a surjective (rational) holomorphic map with connected fibres.

\begin{conjecture}\label{cc+u} 1. ($`C_{n,m}$ conjecture' of Iitaka) For any fibration $f:X\to Y$, $\kappa(X)\geq\kappa(X_y)+\kappa(Y)$. In particular, $\kappa(X)\geq\kappa(Y)$ if $\kappa(X)\geq 0$, since then $\kappa(X_y)\geq 0$. Here $X_y$ is the generic fibre of $f$.

2. (`Uniruledness conjecture'). If $\kappa(X)=-\infty$, $X$ is uniruled (the converse is easy).
\end{conjecture}

{\bf Sketch of proof:} Assume first that $\kappa(X)\geq 0$. The $C_{n,m}$ conjecture, applied to any $f:X\to Y$, directly implies the result in this case.
In general, let $r:X\to R$ be the `rational quotient'. If $X$ is rationally connected, $R$ is a point, and $\kappa_+(X)=-\infty$, since any $Y\leq X$ is uniruled. Thus the equality. Otherwise, let $f:X\to Y$ be any rational fibration. If the generic fibre $X_r$ does not map to a point, $Y$ is uniruled, and $\kappa(Y)=-\infty$. Thus $Y\leq R$ if $\kappa(Y)\geq 0$, which is the claim $\square$

\

\begin{remark} The invariant $\kappa_+(X)$ is  `external' in the sense that it uses manifolds $Y$ others than $X$. We shall now introduce a second, closely related, invariant, which is `internal' to $X$, because it refers only to data defined on $X$ itself.\end{remark}

\

 Let, if $f:X\to Y_p, p:=dim(Y)>0$ is a `fibration', $L_f$ be the line bundle on $X$ defined by: $L_f:=f^*(K_Y)\subset \Omega_X^p$. Thus, $\kappa(X,L_f)=\kappa(Y)$, and $m.L_f\subset Sym^m(\Omega_X^p), \forall m>0$. We may `saturate' $m.L_f$ in $\Sym^m(\Omega^p_X)$ the space of sections of $m.L_f$, and define, for any $L\subset \Omega^p_X$:

\begin{definition}\label{defkL}  Let $L\subset \Omega_X^p$ a rank-$1$ coherent subsheaf, and for any $m>0$, let $H^0(X,m.L)\subset \bar H^0(X,m.L)\subset H^0(X, Sym^m(\Omega_X^p))$ be the subspace of sections taking values in $m.L\subset Sym^m(\Omega_X^p)$ at the generic point of $X$. Thus: $\bar H^0(X,m.L)$ is also the space of sections of the saturation of $m.L$ in $\Sym^m(\Omega^p_X)$.

Thus: $H^0(X,m.L)\subset \bar H^0(X,m.L)$. Let $\bar h^0:=dim_{\Bbb C}\bar H^0$, and :$$\kappa^*(X,L):=\limsup_{m>0}\Big\{\frac{log (\bar h^0(X,m.L))}{log(m)}\Big\}$$
\end{definition}

By standard arguments, one shows that: $\kappa^*(X,L)$ is either $-\infty$ or an integer at most $n$. A fundamental theorem of Bogomolov (see \cite{Bo78}) actually asserts that $\kappa^*(X)\leq p$, with equality iff $L=L_f$ for some dominant rational fibration $f:X\dasharrow Y_p$. However, $Y$ does not need to be of general type in this situation, since due to taking saturation, $\kappa^*(X,L_f)\geq \kappa(X,f^*(K_Y))=\kappa(Y)$, the first inequality being strict in many cases. The difference will be geometrically described below.

\

\begin{definition}\label{defk++} For any $X$, let $\kappa_{++}(X):=max_{\{p>0,L\subset \Omega_X^p, rk L=1\}}\{\kappa^*(X,L)\}$.
\end{definition}

$\kappa_{++}(X)$ is a birational invariant, with: $n\geq \kappa_{++}(X)\geq \kappa_+(X)\geq \kappa(X)$.

\begin{conjecture}\label{ck++} For any $X$, $\kappa_{++}(X)=\kappa_+(X)=\kappa(R(X))$.
\end{conjecture}

When $X$ is rationally connected, it is easy to see that $\kappa_{++}(X)=-\infty$ by restricting $\Omega^p_X$ to a rational curve $C$ with ample normal bundle $N$, and considering its natural filtration with quotients $(\Omega^1_C)\otimes (N^*)^{\otimes k}$. A relative version of this permits to show that $\kappa_{++}(X)=\kappa_{++}(R(X))$ for any $X$. See Lemma \ref{ldesc} below. One is thus reduced, by conjecture \ref{cc+u}(2), to the special case where $\kappa(X)\geq 0$.

 Notice that here, however, the case $\kappa(X)\geq 0$ cannot be derived from the $C_{n,m}$ conjecture, since a geometric interpretation of $\kappa_{++}(X)$ is lacking. Working in a larger category will permit at the same time to give a geometric interpretation of $\kappa^*(X,L_f)$, to formulate a suitable version of the $C_{n,m}$ conjecture, and to give a canonical birational decomposition of any $X$ in terms of `pure' manifolds, for which the canonical bundle has one of  the three basic possible `signs' ($+,0,-$), in some suitable birational sense.

 \section{\bf Extension to the category of `smooth orbifolds'.}

 Let $f:X\to Y_p$ be a fibration, and $L_f=f^*(K_Y)\subset \Omega_X^p$.  We shall always assume that $f$ is `neat' (i.e: that the discriminant locus of $f$ is of snc, and that the $f$-exceptional divisors of $X$ are also $u$-exceptional for some birational map $u:X\to X'$, with $X'$ smooth). This condition can always be realised, by means of Raynaud flattening theorem, after suitable modifications of $X$ and $Y$.

 \

 The invariant $\kappa^*(X,L_f)\geq \kappa(X,L_f)=\kappa(Y)$ can be interpreted geometrically as follows in terms of the `base orbifold' of $f$.
 
 \
  
 A lc pair $(X\vert \Delta)$ consisting of a projective manifold $X$ and of an effective $\bQ$-divisor $\Delta=\sum_{j\in J}a_j.D_j$ , with $a_j=(1-\frac{1}{m_j})$ will be said to be `smooth' if $Supp(\Delta)=\lceil\Delta\rceil$ is of snc (i.e: of simple normal crossings). We shall write $a_j=(1-\frac{1}{m_j})$, or equivalently: $m_j:=(1-a_j)^{-1}\in \Bbb Q\cap\{+\infty\}$, for the $\Delta$-`multiplicity' of $D_j$ (equal to $1$ if $D$ is not one of the $D_j's$). We shall also call such a pair a `smooth orbifold'. They interpolate between the `compact' case in which $\Delta=0$, and the `open' or `logarithmic' case, in which $\Delta=\lceil \Delta\rceil\neq 0$.

 \begin{definition}\label{dbo} The `base orbifold' of $f:X\to Y$ is the pair $(Y\vert \Delta_f)$, with $\Delta_f:=\sum_{E}(1-\frac{1}{m(f,E)}).E$, $E$ running through the set of all prime divisors of $Y$. We define $m(f,E):=inf_{k\in K(E)}\{t_{k,E}\}$, and $t_{k,E}$ by the equality:  $f^*(E)=\sum_{k\in K(E)}t_{k,E}.D_k+R$, $K(E)$ being the set of prime divisors $D_k\subset X$ such that $f(D_k)=E$, while $R$ is $f$-exceptional.
 
 Notice that the sum defining $\Delta_f$ is, in fact, finite, since $m(f,E)=1$ when $E$ is not a component of the discriminant locus of $f$.
 \end{definition}

 The pair $(Y\vert \Delta_f)$ is thought of as a `virtual ramified cover' of $Y$ eliminating by base change the multiple fibres of $f$ in codimension one.  
 
 \

 The geometric interpretation of $\kappa^*(X,L_f)$ is now the following:

 \begin{theorem}\label{tk} ([\cite{Ca07}) For $f$ as above, one has: $\kappa^*(X,L_f)=\kappa(Y,K_Y+\Delta_f)=:\kappa(Y\vert \Delta_f)$.
 \end{theorem}

 The origin of the difference $\kappa^*(X,L_f)-\kappa(X,L_f)$ thus lies in the multiple fibres of $f$. This theorem completes some of the results of \cite{Bo78}. The study of the invariant $\kappa_{++}(X)$ thus leads to the consideration of `smooth' pairs $(X\vert \Delta)$, but for reasons different from the ones in the LMMP.
 
 \

 These `smooth pairs' can be naturally equipped with lots of geometric invariants not considered in the LMMP. We shall briefly list, but not define them\footnote{See \cite{Ca07} for the definitions.}:
 
 \

 $\bullet$ {\it morphisms and birational maps.} We thus obtain a category. If $V=X-D$ is a quasi-projective manifold, with smooth compactification $X$ and complement $D$ such that $(X\vert D)$ is smooth, then the birational class of $V$ does not depend on the compactifying pair $(X\vert D)$ in this category.
 
 $\bullet$ {\it Sheaves of symmetric differentials.} These are locally free sheaves $S^m(\Omega^p(X\vert\Delta))$ interpolating between $Sym^m(\Omega^p_X)$ and $Sym^m(\Omega^p_X(Log(Supp(\Delta))$. When $p=1$, in local analytic coordinates $(x_1,\dots,x_n)$ `adapted' to $\Delta$ (i.e: in which the support of $\Delta$ is contained in the union of coordinate hyperplanes, the hyperplane $x_j=0$ having coefficient $0\leq a_j\leq 1)$, the sheaf $S^m(\Omega^1(X\vert\Delta))$ is generated, as an $\sO_X-$module, by the elements $dx^{(N)}:=\otimes_{j=1}^{j=n}\frac{dx_j^{\otimes N_j}}{x_j^{[a_j.N_j]}}$, parametrised by the $n$-uplets $(N)=(N_1,\dots,N_n)$ such that $m=N_1+\dots+N_n$. 
 
 In particular, $m.(K_X+\Delta)=S^m(\Omega^n(X\vert \Delta))$.

 Morphisms $f:(X\vert \Delta)\to (Y\vert \Delta_Y)$ functorially induce maps  of sheaves of symmetric differentials, moreover, the spaces $H^0(X,S^m(\Omega^p(X\vert\Delta))$) are birational invariants of the smooth pair $(X\vert \Delta)$.

 \
 
 $\bullet$ {\it The `integral' case.} When $\Delta$ is moreover `integral' (i.e: if all $a_j's$ are of the 'standard' form $a_j=1-\frac{1}{m_j}$ with $m_j$ either integral or $+\infty$, that is: $m_j=1$), one can define additionally the $3$ following invariants:  $\pi_1(X\vert \Delta)$, the Kobayashi pseudometric $d_{(X\vert \Delta)}$ on $X$, and the notion of integral points (over any field of definition).
 
 \
 
 The invariants defined above permit to extend, as follows, to `smooth pairs' $(X\vert \Delta)$ the birational invariants $n\geq \kappa_{++}(X\vert \Delta)\geq \kappa_{+}(X\vert \Delta)\geq \kappa_{}(X\vert \Delta)$.

 \

 Let, indeed, a smooth pair $(X\vert \Delta)$ be given.

 \

 $\bullet$ For any $L\subset \Omega^p_X$, and $m>0$, let $ \bar H^0(X\vert \Delta,m.L)\subset H^0(X, S^m(\Omega^p(X\vert \Delta)))$ be the subspace of sections taking values in $m.L\subset S^m(\Omega^p(X\vert \Delta))$ at the generic point of $X$. Equivalently,   $\bar H^0(X\vert \Delta,m.L)= H^0(X, \overline{m.L}^{\Delta})$ is the space of sections of the saturation $\overline{m.L}^{\Delta}$ of $m.L$ in the sheaf $S^m(\Omega^p(X\vert \Delta))$.

 $\bullet$ Define next: $$\kappa^*(X\vert \Delta,L):=\limsup_{m>0}\Big\{\frac{log (\bar h^0(X\vert \Delta,m.L))}{log(m)}\Big\}.$$

$\bullet$  For any neat fibration $f:X\to Y$, we define an orbifold base $(Y\vert \Delta_{(f\vert \Delta)})$ exactly as above when $\Delta=0$, simply replacing there $m(f,E)$ by $m((f\vert \Delta),E):=inf_{k\in K(E)}\{t_{k,E}.m_{\Delta}(D_k)\}$, recalling that $m_{\Delta}(D_k)$ is the multiplicity of $D_k$ in $\Delta$. The notations are those of definition \ref{dbo} above.  The reason for this definition comes from a formula to compute the orbifold base of the composition of two fibrations. When $f$ is only rational, replace it first by a `neat' model.

 The above theorem \ref{tk} still holds: $\kappa^*(X\vert \Delta,L_f)=\kappa(Y\vert \Delta_{(f\vert \Delta)})$.
 
 \

$\bullet$ Define finally: $\kappa_+(X\vert \Delta):=max_{\{Y\leq X\}}\{\kappa(Y\vert \Delta_{(f\vert \Delta)}))\},$ and:

                                      \

                                      $\kappa_{++}(X\vert \Delta):=max_{\{p>0,L\subset \Omega_X^p, rk L=1\}}\{\kappa^*(X\vert \Delta,L)\}$
.

 \

 The conjecture \ref{ck++} can now be partially extended to `smooth orbifolds':
 
 \begin{conjecture}\label{ck++d} For any `smooth pair' $(X\vert \Delta)$ such that $\kappa(X\vert \Delta)\geq 0$, one has:  $\kappa_{++}(X\vert \Delta)=\kappa(X\vert \Delta)$.
  \end{conjecture} 
 
In general, we shall conjecture that $\kappa_{++}(X\vert \Delta)=\kappa_+(X\vert \Delta)=\kappa(R^*\vert \Delta_{(r^*\vert\Delta)})$, the fibration $r^*:(X\vert \Delta)\to R^*$, which is a substitute of the rational quotient, being conditionally defined in proposition \ref{prq} when $\kappa(X\vert \Delta)=-\infty$.

\

We shall provide in \S \ref{oratc} a conjectural geometric interpretation (see \ref{corc} below) of the conditions $\kappa=-\infty$ and $\kappa_{++}=-\infty$ in the orbifold context, in terms of `orbifold' rational curves.

\

\section{\bf Orbifold Additivity.}

Let $(X\vert \Delta)$ be a `smooth' pair, and $f:X\to Y$ be a `neat' fibration, with `orbifold base' $(Y\vert \Delta_{(f\vert \Delta)})$. Notice that the restriction of $\Delta$ to a generic fibre $X_y$ of $f$ induces a smooth pair $(X_y\vert \Delta_y)$.

\begin{conjecture} \label{cc} (`$C_{n,m}^{orb}$-conjecture') $\kappa(X\vert \Delta)\geq \kappa(X_y\vert \Delta_y)+\kappa(Y\vert \Delta_{(f\vert\Delta)})$.
\end{conjecture}

Remark that, even when $\Delta=0$, this strengthens Iitaka Conjecture $C_{n,m}$ (because of the second term on the right hand side, which takes multiple fibres into account).

\begin{theorem}\label{tgt} (\cite{Ca07}) When the orbifold base of $(f\vert \Delta)$ is of general type (ie: if $\kappa(Y\vert \Delta_{(f\vert\Delta)})=dim(Y))$, we have: $\kappa(X\vert \Delta)\geq \kappa(X_y\vert \Delta_y)+dim(Y)$.
\end{theorem}

The proof is an orbifold adaptation of Viehweg's arguments used when $\Delta=0$. Nevertheless, the orbifold context considerably extends the range of applicability. Applications wil be given in \S \ref{sc+dec}. We first derive some (conditional) conclusions of the conjecture.

We now introduce the two fundamental fibrations of birational classification in the orbifold context.

The first one is the Iitaka-Moishezon fibration $J:(X\vert \Delta)\to J(X\vert \Delta)$, defined by a suitable linear system $m.(K_X+\Delta)$ when $\kappa(X\vert \Delta)\geq 0$. Its two defining properties are that its generic orbifold fibres $(X_j\vert \Delta_j)$ have $\kappa=0$, and that its base dimension is $\kappa(X\vert \Delta)\geq 0$. 

\

Applying $C_{n,m}^{orb}$ to the orbifold Iitaka-Moishezon fibration gives a partial answer to conjecture \ref{ck++d}:

\begin{proposition}\label{pk=o} Assume that the $C_{n,m}^{orb}$ conjecture \ref{cc} holds, then $\kappa^{*}(X\vert \Delta, L_f)\leq\kappa(X\vert \Delta)$ for any fibration $f:X\dasharrow Y$.
\end{proposition}

The second fundamental fibration is a weak (conditional) version of the `rational quotient'. Its existence requires assuming $C_{n,m}^{orb}$.

\begin{proposition}\label{prq}(\cite{Ca07})  Assume $C_{n,m}^{orb}$. For any smooth $(X\vert \Delta)$, there exists a (birationally) unique fibration $r^*:(X\vert \Delta)\to R^*:=R^*(X\vert \Delta)$ such that:

1. Its generic orbifold fibres have $\kappa_{+}=-\infty$.

2. Its orbifold base has $\kappa\geq 0$, or is a point iff $\kappa_{+}(X\vert \Delta)=-\infty$.
\end{proposition}

We now reformulate the general (conditional) version of conjecture \ref{ck++d}, in complete analogy with the case $\Delta=0$:

\begin{conjecture}\label{ck++d'} Assume $C_{n,m}^{orb}$ (it is needed to define $r^*$). For any smooth $(X\vert \Delta)$, one has: $\kappa_{++}(X\vert \Delta)=\kappa(R^*\vert \Delta_{(r^*\vert \Delta)})$. Here $(R^*\vert \Delta_{(r^*\vert \Delta)})$ is simply the orbifold base of (any neat model of) $r^*:(X\vert \Delta)\dasharrow R^*$.
\end{conjecture}

\begin{remark} Although the fibration $r^*$ and $R^*$ are well-defined up to birational equivalence, it is not known whether its orbifold base $(R^*\vert \Delta_{(r^*\vert \Delta)})$ is uniquely defined up to birational equivalence. Its Kodaira dimension is however well-defined, independently of the choices made. See \cite{Ca07}.
\end{remark}

We shall now reduce this fibration to a composition of fibrations of the LMMP, and the conjecture \ref{ck++} to  some more standard conjectures.

\section{\bf Reduction to $2$ other conjectures.}

\

We now formulate three other conjectures, the first two ones being standard in the LMMP (due to V. Shokurov and C. Birkar).

\begin{conjecture}\label{c5} Let $(X,\Delta)$ be a l.c pair. 

1. There exists a sequence of divisorial contractions and flips $s:X\dasharrow X'$ such that if $\Delta'=s_*(\Delta)$, then either $K_{X'}+\Delta'$ is nef, or there exists a fibration $f:X'\to Y'$ with Fano, positive-dimensional orbifold fibres $(X'_y,\Delta'_y)$.

1'. If $K_{X'}+\Delta'$ defined above is nef, it is $\Bbb Q$-effective. (It is a weak form of the Abundance conjecture, formulated in \ref{abcj} below).

2. If $c_1(X\vert \Delta)=0$, if $m>0$ is an integer, and if $C=H_1\cap\dots H_{n-1}$ is a general Mehta-Ramanathan curve on $X$, the restriction of $\otimes_{h}S^{m_h}(\Omega^{p_h}(X\vert \Delta))$ to $C$ is semi-stable (i.e: all of its subsheaves have nonpositive degree), for any finite sequence $(m_h,p_h)$ of pairs of positive integers .
\end{conjecture}

The first conjecture is known for klt pairs, by [BCHM], the second one a special case of the Abundance conjecture, the third one is simply the orbifold version of Miyaoka's generic semi-positivity theorem.

Let us give first a description of the fibration $r^*$ using the conjectures \ref{c5}(1,2).

\begin{definition}\label{dr} Let $(X\vert \Delta)$ be a smooth projective orbifold. Define the fibration $r=r_{(X\vert \Delta)}:(X\vert \Delta)\dasharrow (Y\vert \Delta_Y)$ to be:

1. The (orbifold) identity map if $\kappa(X\vert \Delta)\geq 0$. 

2. Any neat model of the composition map $r:=(f'\circ s):X\dasharrow Y'$ of \ref{c5}(1) if $\kappa(X\vert \Delta)=-\infty$, with orbifold base $(Y\vert \Delta_Y)$.\end{definition}

Notice that neither $r$, nor $(Y\vert \Delta_Y)$ are uniquely defined, up to birational (orbifold) equivalence. Nevertheless, the composition $r^n, n=dim_{\Bbb C}(X),$ is well-defined, by the following:

\begin{theorem}\label{tr^*} Let $(X\vert \Delta)$ be a smooth $n$-dimensional projective orbifold. Assume that the conjectures \ref{c5} and $C_{n,m}^{orb}$ hold. Then $r^*=r^n$ (for any possible choice of the sequence of $r's)$.
\end{theorem}

{\bf Proof:} Because of the uniqueness of the map $r^*$ (up to birational equivalence), we simply need to show that, on any neat model of $r^n$, we have: $\kappa_+=-\infty$ for its general orbifold fibre, and $\kappa\geq 0$ for its orbifold base (on some neat model). If $\kappa(X\vert \Delta)\geq 0$, $r^n$ is the identity map, and the claim is obvious (assigning $\kappa_+=-\infty$ to orbifold points). Otherwise, let $m>0$ be the smallest integer such that $\kappa(Y\vert \Delta_Y)\geq 0$, for the base orbifold of some sequence of $r's$ of lentgh m, and of composition $r^m$. We then have $m\leq n$, since the dimension decreases at each of the $m$ steps, since the intermediate orbifold bases have $\kappa=-\infty$. Now, the intermediate fibrations have generic fibres which are birational to l.c Fano orbifold pairs, and have thus $\kappa_+=-\infty$, by lemma \ref{fanok++}, proved below. Since a composition of rational fibrations with general orbifold fibres having $\kappa_+=-\infty$ also has this property (by \cite{Ca07}, 7.14), we conclude that the general orbifold fibres of $r^m$ also have $\kappa_+=-\infty$ $\square$

\

A second consequence of the conjecture \ref{c5} is:

\begin{theorem}\label{trlmm} (\cite{Ca 09}, theorem 10.5) Assume $C_{n,m}^{orb}$ and conjecture \ref{c5}. Then conjecture \ref{ck++d} also holds.
\end{theorem}

A stronger (numerical dimension) version will be proved in detail in \ref{tnu++} below, along parallel lines. We thus do not give the proof again here.

\section{\bf The core and its canonical decomposition.}\label{sc+dec}

\

\begin{definition}\label{ds} We say that $(X\vert \Delta)$ is `special' if $\kappa^*(X\vert \Delta, L)<p$ for any $p>0$ and any rank-$1$ coherent subsheaf $L\subset \Omega^p_X$. 

Equivalently, this means that there is no fibration $f: (X\vert \Delta)\dasharrow Y$ such that its orbifold base is of general type (on any holomorphic neat model), with $dim(Y)>0$.
\end{definition}

\

\begin{remark} Rank-$1$ coherent subsheaf $L\subset \Omega^p_X, p>0$ of maximum positivity (i.e: with $\kappa^*(X\vert \Delta, L)=p$) are called `$\Delta$-Bogomolov sheaves'. Being special thus means that there are no such sheaves on $(X\vert \Delta)$.
\end{remark}

\

Special orbifolds are natural higher-dimensional generalisations of rational and elliptic curves, with the same expected qualitative properties. They are the exact opposite of orbifolds of  `general type'.

\

\begin{corollary}(of theorem \ref{tgt}) If $\kappa(X\vert \Delta)=0$, or if $(X\vert \Delta)$ is Fano (i.e: if $-(K_X+\Delta)$ is ample, then $(X\vert \Delta)$ is special.
\end{corollary}

By the very definition, $(X\vert \Delta)$ is special if $\kappa_{++}(X\vert \Delta)=-\infty$. It is unknown whether $\kappa_{++}(X\vert \Delta)=-\infty$ if $(X\vert \Delta)$ is Fano. This follows however from conjecture \ref{c5}(2), as we have seen above.

\

The next result describes unconditionally the structure of arbitrary smooth orbifolds, in terms of its antithetic maximal parts (special `subobjects' vs general type `quotients')
:

\begin{theorem}\label{tc} For any smooth pair $(X\vert \Delta)$, there is a (birationally) unique functorial fibration $c:(X\vert \Delta)\to C(X\vert \Delta)=C$ such that:

1. Its general fibres are special.

2. Its orbifold base is of general type (or a point, iff $(X\vert \Delta)$ is special).

This fibration is called `the core' of $(X\vert \Delta)$. \end{theorem}

\begin{remark} The fibration $c$ is determined by the (unique) $\Delta$-Bogomolov sheaf $L\subset \Omega^p_X$ on $(X\vert \Delta)$, with $p>0$ maximum.

\end{remark}

The second (conditional in $C_{n,m}^{orb}$) structure theorem is:

\begin{theorem}\label{tc=jr} Assume the preceding proposition \ref{prq} (since it uses $C_{n,m}^{orb}$). Then, for any smooth pair $(X\vert \Delta)$, the core map of $(X\vert \Delta)$ has the following decomposition as a composition  of $2n$ canonically defined fibrations: $c=(J\circ r^*)^n$.

In particular, $(X\vert \Delta)$ is special iff it is a tower of fibrations with general orbifold fibres having either $\kappa_{+}=-\infty$, or $\kappa=0$.
\end{theorem}

Notice that, even if we are only interested in the case $\Delta=0$, non-trivial orbifold divisors will usually appear in the above decomposition.

This (conditional) decomposition often permits to reduce the study of arbitrary  manifolds to that of smooth pairs of the three basic types: $\kappa_{+}=-\infty$, $\kappa=0$, or of general type. It naturally leads to a conjectural extension of S. Lang's conjectures in arithmetics and complex hyperbolicity, for all manifolds and even smooth orbifolds. See [Ca 01] and [Ca 07] for details.

\section{\bf Numerical dimension version.}

\

Let, in this section, $X$ be a complex projective connected, $n$-dimensional $\Bbb Q$-factorial normal space, $A$ and $D$ be $\Bbb Q$-divisors on $X$, with $A$ ample. 

\

The {\it numerical dimension of $D$} is defined as the following real number: $\nu(X,D):=sup_{k\geq 0}\Big\{limsup_{m>0}\Big\{\frac{Log(h^0(X,mD+A))}{Log(m)}\Big\}\Big\},$ for $m>0$ integral and sufficiently divisible. 

Easy standard arguments show that:

1. $\nu(X,D)=-\infty$, or is real, and lies in $[ 0,n]$. 

2. $\nu(X,D)$ does not depend on the choice of $A$. 

3. $\nu(X,D)\geq \kappa(X,D)$.

4. $\nu(X,D)=-\infty$ if and only $D$ is not pseudo-effective (this is one of the definitions of pseudo-effectivity).

5. When $D$ is nef, it is an easy consequence of Kodaira vanishing and Riemann-Roch that $\nu(X,D)=limsup_{m>0}\Big\{\frac{Log(h^0(X,mD+A))}{Log(m)}\Big\}=\nu'(X,D)\in \{0,1,\dots,n\}$, for any ample $A=K_X+(n+2)H$, where $H$ is any ample line bundle on $X$, and where $\nu'(X,D)$ is the largest integer $\ell$ such that $D^{\ell}\in H^{2.\ell}(X, \bZ)$ is not numerically zero. The Kodaira vanishing indeed says that $h^0(X,mD+A)=\chi(X, \sO_X(mD+A))$, for any $m\geq 0$. When $D$ is only assumed to be pseudo-effective, the Nadel vanishing theorem implies the same equality, but only after tensorising $\sO_X(mD+A)$ with the multiplier ideal sheaf $\sI(mD+A)$, which cannot be controlled without further ideas. 

6. One may however wonder whether $\nu(X,D)$ is not an integer, if nonnegative, and if $\nu(X,D)=\nu_A(X,D)$ for $A$ sufficiently ample (for example $A=K_X+(n+2).H$, as above). And also what is the relationship between $\nu(X,D)$ and the numerical dimension of $D$ defined by N. Nakayama in \cite{N} and S. Boucksom in \cite{Bo}.

\

One form of the so-called `Abundance Conjecture' is the following:

\begin{conjecture}\label{abcj} Assume $(X\vert D)$ is a `log-canonical pair'. Then:

 $\nu(X,K_X+\Delta)=\kappa(X,K_X+\Delta)$.
\end{conjecture}

This is known when $D$ is `big' and $(X\vert D)$ is klt (\cite{BCHM} and \cite{Pa}). This is also known when $\nu(X,K_X+\Delta)=0$ if $q(X)=0$, as follows from \cite{N} and \cite{Bo}. When $\Delta=0$, the case $q\geq 0$ follows from a more general statement in \cite[\S 3]{CP05}. Finally, the general lc case when $\nu=0$ is established (in a more general form) in \cite{CKP}, using the purely logarithmic case proved in \cite{Kaw10}.

\begin{proposition}\label{pabiad} Assume conjecture \ref{abcj}. Then Conjecture $C_{n,m}^{orb}$ is true.
\end{proposition}

{\bf Proof:} See \cite{Ca 09}, \S 10 for a proof using the weak positivity of the direct images of the orbifold pluricanonical sheaves. We give in \S 8 below a simple proof in the particular case where $\Delta=0$, using the pseudo-effectivity of $f_*(K_{X/Y})$ when $X_y$ is not uniruled $\square$

\

We shall now state and conditionally prove a `numerical dimension' version of theorem \ref{trlmm}. For this we first need to define the `numerical dimension' version of $\kappa_{++}$.

\begin{definition}\label{dnu} Let $E_.=(E_m)_{m\in \Bbb N^{>0}}$ be a family of vector bundles on $X$, together with generically isomorphic bundle maps $S^m(E_1)\to E_{m}$ for any $0<m \in \Bbb N$. Let $L\subset E_1$ be a rank-one coherent susheaf, and let $\overline{m.L}^{E_.}$ be the saturation of the image of $Sym^m(L)$ in $E_m$, for any $m>0$.

Let $A$ be an ample line bundle on $X$. 

We define $\nu_A(X\vert E_.,L):=limsup_{m>0}\Big\{\frac{Log(h^0(X,\overline{m.L}^{E_.}\otimes A))}{Log(m)}\Big\}$, and:

 $\nu(X\vert E_.,L):=max_{k>0}\{\nu_{kA}(X\vert E_.,L)\}$.
\end{definition}

\

Of course, we always have: 

1. $\nu(X,D)=-\infty$, or is real, and lies in $[ 0,n]$. Indeed $\nu(X\vert E_.,L)$ is bounded by the maximum dimension of the image of $X$ by the rational  map deduced from any nonzero linear system $h^0(X,\overline{m.L}^{E_.})$.

2. $\nu(X\vert E_.,L)$ does not depend on the choice of $A$.

3. $\nu(X\vert E_.,L)\geq \nu_A(X\vert E_.,L)\geq \kappa(X\vert E_.,L):=limsup_{m>0}\Big\{\frac{Log(h^0(X,\overline{m.L}^{E_.}))}{Log(m)}\Big\}$

\

The main examples considered here are:

\begin{example} 1. $L=E_1$, and $E_m:=m.E_1$. This is the standard case.

2. Let $(X\vert \Delta)$ be a smooth orbifold, $p>0$, and $E_m:=S^m(\Omega^p(X\vert \Delta))$. In this case, $\nu_A(X\vert E_.,L)$ is denoted by: $\nu_A(X\vert \Delta, L)$, and similarly for $\nu(X\vert \Delta, L)$. We also denote $\overline{m.L}^{E_.}$ by $\overline{m.L}^{\Delta}$ in this case.
\end{example}

\

\begin{definition} If $(X\vert \Delta)$ is a smooth orbifold, then we define: 

$\nu_{++}(X\vert \Delta)=max_{\{p>0,L\subset \Omega^p_X\}}\nu(X\vert \Delta, L)$.
\end{definition}

\

We now have the following strengthening of theorem \ref{trlmm}:

\

\begin{theorem}\label{tnu++} Assume conjectures \ref{abcj} and \ref{c5}. Then, for any smooth orbifold $(X\vert \Delta)$, one has: $\nu_{++}(X\vert \Delta)= \kappa_{++}(X\vert \Delta)=\kappa(R^*(X\vert \Delta)\vert \Delta_{(r^*\vert \Delta)})$.
\end{theorem}

Let us remark that, although the base orbifold $(R^*(X\vert \Delta)\vert \Delta_{(r^*\vert \Delta)})$ is not known to be birationally wel-defined, its canonical dimension $\kappa$ is well-defined.

\

{\bf Proof:}  Let $(X\vert \Delta)$ be a smooth orbifold with $X$ projective. Let $\sF\subset \Omega^p_X$ be a rank-one coherent subsheaf.

Combining conjectures \ref{c5} and \ref{abcj}, we first observe that theorem \ref{tnu++} holds when $\kappa(X\vert \Delta)=0$. Indeed: we can assume, using the birational map $s:(X\vert \Delta)\dasharrow (X'\vert \Delta')$ provided by conjecture \ref{c5}(1), with $c_1(X'\vert \Delta')=0$, in which case the claim immediately follows from conjecture \ref{c5}(2) by restricting to a general Mehta-Ramanathan curve $C\subset X'$, by means of the following lemma \ref{lrc} (1). From the following lemma \ref{lrc} (2), we now deduce theorem \ref{tnu++}  also when $\kappa(X\vert \Delta)\geq 0$, by using a neat model of the Moishezon-Iitaka fibration for $(X\vert \Delta)$.

\begin{lemma}\label{lrc} Let $X$ be smooth projective connected, and $E_.$ and $L\subset E_1$ be as above. Let $C=H_1\cap\dots\cap H_{n-1}\subset X$ be a general Mehta-Ramanathan curve on $X$ of genus $g(C)$.

1. Assume that $L.C\leq 0$, and that $H_i.C>A.C$, for $i=1,\dots,(n-1)$.
Then, for each ample $A$ on $X$, and for any set of $(2.g(C)+(A.C))$ distinct points $c_k$ on $C$, the natural restriction map: 
$H^0(X,\overline{m.L}^{E_.}\otimes A)\to \oplus_{k=1}^{k=2.g(C)+A.C} (m.L+A)_{\vert c_k}$ is injective. In particular, $h^0(X,\overline{m.L}^{E_.}\otimes A)\leq (2.g(C)+A.C)$ for any $m>0$, and $\nu_A (X\vert E_.,L)\leq 0$.

2. If $f:X\to Y$ is a fibration such that, for any integer $k>0$, there exists a bound $B(k)$ such that $h^0(X_y, \overline{(m.L+k.A)}^{ E_.}_{\vert X_y})\leq B(k)$ for any $m>0$, then $\nu (X\vert E_.,L)\leq p:=dim(Y)$.
\end{lemma}

{\bf Proof:} 1. It is sufficient to show that the kernel $Ker$ of the restriction map $res: H^0(X,\overline{m.L}^{E_.}\otimes A)\to H^0(C,\overline{m.L}^{E_.}\otimes A)_{\vert C}$ is zero, since the evaluation map on the points $c_k$ is injective. But $Ker=H^0(X,\overline{m.L}^{E_.}\otimes A\otimes \sI_C)$, where $\sI_C\cong \oplus_{i=1}^{i=(n-1)} \sO_X(-H_i)$ is the ideal of $C$. Thus $Ker=\{0\}$, since $(m.L+A-H_i).C<0$ for all $i's$, and $C$ belongs to an $X$-covering family of curves of $X$.

2. It is sufficient to show that $\nu_A (X\vert E_.,L)\leq p$, and then to replace $A$ by $k.A$ in the argument. Let $Z:=H_1\cap\dots\cap H_{n-p}$ be the smooth connected complete intersection of very ample divisors on $X$, such that the degree of the restricted map $f_{\vert Z}:Z\to Y$ is at least $B(1)+1$, and $Z\cap X_y:=Z_y$ consists of a $B(1)+d,d>0$ points $c_{k,y}$ in general position on $X_y$. The restriction map: $H^0(X_y,\overline{m.L}^{E_.}\otimes A)_{X_y})\to \oplus_{k=1}^{k=B(1)+d}(m.L+A)_{c_{k,y}}$ is thus injective, and so is thus the restriction map $H^0(X,\overline{m.L}^{E_.}\otimes A)\to H^0(Z,\overline{m.L}^{E_.}\otimes A)_{\vert Z})$. Thus $\nu_A (X\vert E_.,L)\leq \nu_A (Z\vert (E_.)_{\vert Z},L_Z)\leq p=dim(Z)$ $\square$

\

The general case will result from the following:

\begin{proposition}\label{pdesc} Let $(X\vert \Delta)$ be smooth. Then $\nu_{++}(X\vert \Delta)=\nu_{++}(Y\vert \Delta_Y)$, if $(Y\vert \Delta_Y)$ is the orbifold base of any neat representative of $r^*:(X\vert \Delta)\to R^*(X\vert \Delta)$.
\end{proposition}

By the preceding arguments, and assuming conjecture \ref{abcj}, this proposition indeed implies that $\nu_{++}(X\vert \Delta)=\nu_{++}(Y\vert\Delta_Y)=\kappa(Y\vert \Delta_Y)=\kappa(R^*(X\vert \Delta)\vert \Delta_{(r^*\vert \Delta)})$, which is what theorem \ref{tnu++} claims, since $\kappa(Y\vert\Delta_Y)\geq 0$, for $(Y\vert \Delta_Y)$ as in \ref{pdesc}. 

\

{\bf Proof (of \ref{pdesc}):} Since, by theorem \ref{tr^*}, we have $r^*=r^n$, for any length-n composition of rational fibrations $(f'\circ s)$ with log-Fano fibres (in the sense of the statement of conjecture \ref{c5}(1)), it is sufficient to show that the invariant $\nu_{++}$ is preserved under such fibrations.

\

We first establish the statement for smooth pairs $(X\vert \Delta)$ which are birational to Log-Fano pairs.

\begin{lemma}\label{fanok++} Let $g:(X\vert \Delta)\to (X'\vert \Delta')$ be a birational map from the smooth orbifold $(X\vert \Delta)$ to the log-canonical Fano pair $(X'\vert \Delta')$ such that $f_*(\Delta)=\Delta'$. Assume that conjecture \ref{c5}.(2) holds. Then, for any polarisations of $X'$, and any corresponding general Mehta-Ramanathan curve $C\subset X'$, identified with its strict transform on $X$, the following properties hold:

1. For any finite sequence of pairs of nonnegative integers $(N_h,q_h),h=1,\dots, s$, and any rank-one coherent subsheaf $\sF\subset S^{N_1}\Omega^{q_1}(X\vert \Delta)\otimes \dots \otimes S^{N_s}\Omega^{q_s}(X\vert \Delta)$, the restriction $det(\sF)_{\vert C}$ has nonpositive degree at most:  $-[(\sum_{h=1}^{h=s}q_h.N_h)-M.n^2]$, $M$ being any integer such that: $-M.(K_{X'}+\Delta')$ is very ample.

2. $H^0(X,S^{N_1}\Omega^{q_1}(X\vert \Delta)\otimes \dots \otimes S^{N_s}\Omega^{q_s}(X\vert \Delta)\otimes A)=\{0\}$, for any ample line bundle $A$ on $X$,  if $(\sum_{h=1}^{h=s}q_h.N_h)>M.n^2+A.C$.

3. $h^0(X,\overline{m.L}^{\Delta}\otimes A)=0$ if  $m.(\sum_{h=1}^{h=s}q_h.N_h)>M.n^2+A.C$. 

4. $\nu_{++}(X\vert \Delta)=-\infty$.
\end{lemma}

{\bf Proof:} 1. Let $\sG:=S^{N_1}\Omega^{q_1}(X\vert \Delta)\otimes \dots \otimes S^{N_s}\Omega^{q_s}(X\vert \Delta)$.  Let $H'=\sum_{j=1}^{j=n}H_j$, where the $H_j's$ are general members of $M.(-(K_{X'}+\Delta'))$, $M>0$ being a sufficiently large integer, the $H_j's$ being chosen so that $(X'\vert \Delta'+\frac{1}{mn}.H'):=(X'\vert \Delta")$ is log-canonical, with $(K_{X'}+\Delta")$ trivial on $X'$, and such that $\Delta^+:=(\Delta+\frac{1}{mn}.H)$ has normal crossings support
, $H$ being the strict transform of $H'$ in $X$. Choose $C\subset X'$ meeting $H'$ transversally, but not meeting the indeterminacy locus of $g^{-1}$,and so identified with its strict transform on $X$. Then $L\subset \sG^+:=S^{N_1}\Omega^{q_1}(X\vert \Delta^+)\otimes \dots \otimes S^{N_s}\Omega^{q_s}(X\vert \Delta^+)_{\vert C}$ has nonpositive degree, by conjecture \ref{c5} (2), since it is a rank one subsheaf of the locally free sheaf $\sG^+$, assumed to be semi-stable, and with trivial determinant.

Assume now that the $H_j's$ have been chosen in such a way that they build a system of coordinate hyperplanes for suitable local coordinates at a generic point $a\in X'$ outside of the support of $\Delta'$ and the indeterminacy locus of $g^{-1}$, and belonging to the smooth locus of $X'$. We choose also $C$ in such a way that $a\in C$.  The natural inclusion $\sG\subset \sG^+$ now vanishes at order at least $[(\sum_{h=1}^{h=s}q_h.N_h)-M.n^2]$ at $a$, as follows from lemma 3.3 of \cite{Ca 09}, since $q_h\leq n$, for any $h=1,\dots, s$. It follows that the degree of $L$ on $C$ is nonpositive, and is at most: $-[(\sum_{h=1}^{h=s}q_h.N_h)-M.n^2]$, and so that it is negative if $(\sum_{h=1}^{h=s}q_hh.N_h)>M.n^2$.

2. This follows from the fact that the restriction of such a section to $C$ vanishes, unless $N_h.q_h=0, h=1,\dots,s$, since a nonzero section would otherwise generate a (locally free) rank-one coherent sheaf of negative degree on $C$, by the previous estimate on the degree of $\sF_{C}$. The last two assertions are now obvious $\square$

\

We shall now deal with the (rational) fibrations having log-canonical Fano fibres. Let us first describe the situation provided by conjecture \ref{c5}(1), assuming that $\kappa(X\vert \Delta)=-\infty$.
Applying conjecture \ref{c5} (1), we get a birational map $s:(X\vert \Delta)\dasharrow (X'\vert \Delta')$ and a log-Fano fibration $f:(X'\vert \Delta')\to Y$, with $dim(Y)<n$ and $(X'_y\vert \Delta'_y)$ log-canonical and Fano for generic $y\in Y$. We can moreover assume that $s:(X\vert \Delta)\to (X'\vert \Delta')$ is  regular and and a log-resolution, and also that $f\circ s:(X\vert \Delta)\to (Y\vert \Delta_Y)$ is a neat orbifold morphism, by making a suitable orbifold modification of $(X\vert \Delta)$ and choosing appropriate multiplicities on the divisors of $X$ which are $f\circ s$-exceptional.

We are thus in the position to apply the following lemma \ref{ldesc}, which implies the claim, and thus proposition\ref{pdesc} and theorem \ref{tnu++} $\square$

\begin{lemma}\label{ldesc} Let $(X\vert \Delta)$ be a smooth orbifold, and $f:X\to Y$ be a neat fibration which is an orbifold morphism with generic smooth orbifold fibre $(X_y\vert \Delta_y)$ and smooth orbifold base $(Y\vert \Delta_Y:=\Delta_{(f\vert \Delta)})$. 

1. Assume that, for any finite sequence of pairs of positive integers $(N_h,q_h),h=1,\dots, t$, one has: $H^0(X_y,S^{N_1}\Omega^{q_1}(X_y\vert \Delta_y)\otimes \dots \otimes S^{N_t}\Omega^{q_t}(X_y\vert \Delta_y))=\{0\}$. Then $f_*(S^N\Omega^q(X\vert \Delta))=S^N\Omega^q(Y\vert \Delta_{Y})$, for any integer $N>0$ and $q>0$.

2. Assume, additionally, that $\nu_{++}(X_y\vert \Delta_y)=-\infty$. Then, we also have:  $\nu_{++}(X\vert\Delta)=\nu_{++}(Y\vert \Delta_Y)$.
\end{lemma}

{\bf Proof:}  1. This is just lemma 4.23 of \cite{Ca 09}. (The statement given there is global on $X$, but its proof applies locally over $Y$).

2. For any ample line bundle $A$ on $X$, and any pair $(N,q)$ of positive integers, we thus have: $H^0(X,S^N\Omega^q(X\vert \Delta)\otimes A)\cong H^0(Y,S^N\Omega^q(Y\vert \Delta_Y)\otimes f_*(A))$. Assume that some rank-one coherent subsheaf $\sF\subset \Omega^q_X$ is such that $\nu_A(X\vert \Delta,\sF)\geq 0$. 

We shall prove first that there exist $\sG\subset \Omega^q_Y$ such that, generically over $Y$, $\sF=f^*(\sG)$. Otherwise, there would exist a largest $s>0$, such that $\sF$ had nonzero image $\sF_y$ in the quotient $f^*(\Omega^{(q-s)}_Y)\wedge \Omega_{X_y}^{s}\cong ( \Omega_{X_y}^{s})^{\oplus R}, R:=(_p^{q-s})$, of the graduation of the natural filtration of $\Omega^q_{X\vert X_y}$ determined by $f$ on its generic orbifold fibre $(X_y\vert \Delta_y)$ (see \cite{Ca 09}, \S 4).  By the assumption that $\nu_A(X\vert \Delta,\sF)\geq 0$, there are arbitrarily large integers $m$ such that $\overline{(m.\sF)}^{\Delta}\otimes A$ has a nonzero section. But these sections would induce by projection  nonzero sections of $S^m\Omega^{s}(X_y\vert \Delta_y)\otimes A_{X_y}$ contained in $\overline{(m.\sF_y)}^{\Delta_y}$, contradicting the hypothesis that $\nu_{++}(X_y\vert \Delta_y)=-\infty$.

From the preceding arguments, we deduce that, for any $m>0$, we have: $h^0(X, \overline{(m.\sF)}^{\Delta}\otimes A)=h^0(Y, \overline{(m.\sG)}^{\Delta_Y}\otimes f_*(A))$. Let now $B$ be an ample line bundle on $Y$. Since there exists positive integers $k$ and $r$ such that the sheaf $f_*(A)$ can be embedded in $(k.B)^{\oplus r}$, we see that $h^0(Y, \overline{(m.\sG)}^{\Delta_Y}\otimes f_*(A))\leq r. h^0(Y, \overline{(m.\sG)}^{\Delta_Y}\otimes (k.B))$. Thus $\nu_A(X\vert \sF)\leq \nu_{kB}(Y\vert \Delta_Y, \sG)$. Since this holds for any $A$, the lemma is established $\square$

\

\section{\bf Orbifold rational curves.}\label{oratc}

We shall now provide a conjectural geometric interpretation (see \ref{corc} below) of the conditions $\kappa=-\infty$ and $\kappa_{++}=-\infty$ in the orbifold context, in terms of `orbifold' rational curves. This interpretation is entirely similar to the case when $\Delta=0$, once the notion of orbifold rational curves is defined. 

\begin{definition}\label{dorc} ([Ca 07, \S 6]) Let $(X\vert \Delta)$ be a smooth orbifold, with $\Delta:=\sum_{j\in J}(1-\frac{1}{m_j}).D_j$. Let $C$ be a smooth connected projective curve. A map $g:C\to (X\vert \Delta)$ is a $\Delta$-rational (resp. a $\Delta$-elliptic) curve if:

1. It is birational onto its image, which is not contained in $Supp(\Delta)$.

2. $deg(K_{C}+\Delta_g)<0$ (resp. $deg(K_{C}+\Delta_g)=0)$, where $\Delta_g$ is the orbifold divisor on $\Bbb P^1$ which assigns to any $a\in \Bbb P^1$ the multiplicity $1$ if $g(a)\notin Supp(\Delta)$, and otherwise the multiplicity $m_g(a):=max_{j\in J(a)}\{max\{1,\frac{m_j}{t_{j,a}}\}\}$. 

Here $J(a):=\{j\in J\vert g(a)\in D_j\}$, and $t_{j,a}$ is the order of contact of $g$ and $D_j$ at $a$, defined by the equality: $g^*(D_j)=t_{j,a}.\{a\}+\dots$, if $j\in J(a)$.

Notice that $C\cong \Bbb P^1$ if $g$ is $\Delta$-rational, but that $C$ is either rational or elliptic if $g$ is $\Delta$-elliptic. If $C$ is elliptic, it is $\Delta$-elliptic if and only if $g(C)$ does not meet the support of $\Delta$.

There is also a stronger `divisible' version of these notions, which we shall not define here.
\end{definition}

\begin{example}\label{exorc} 1. $\Delta_g=0$, and so $g$ is a $\Delta$-rational curve if $t_{j,a}\geq m_j$ for any $a\in \Bbb P^1$ and $j\in J(a)$. A special case is when $C\subset X$ is a rational curve meeting each of the $D_j$ in distinct smooth points of $Supp(\Delta)$, each with multiplicity $m_j$. In this case, $\Delta_g=0$. Such a rational curve will be said `$\Delta$-nice' in the sequel. In this case, $m_j$ must divide $D_j.C$, for each $j\in J$.

2. If $\Delta=Supp(\Delta)$, so logarithmic, a rational curve on $X$ is $\Delta$-rational (resp. $\Delta$-elliptic)  if and only if its normalisation meets $\Delta$ in at most one point (resp. in exactly two points).

3. If $\Delta'\leq \Delta$, then any $\Delta$-rational (resp. $\Delta$-elliptic) curve is also $\Delta'$-rational (resp. either $\Delta'$-rational or $\Delta'$-elliptic) .
\end{example}

\begin{example}\label{oec} {\bf Orbifold-ramified covers} Let $u:X'\to (X\vert \Delta)$ be a surjective finite ramified cover, with $X'$ smooth connected and $(X\vert \Delta)$ smooth. We say that $u$ is `orbifold-ramified' if the $m_j's$ are integers, if $u$ is unramified over $X-Supp(\Delta)$, and if, for any $j\in J$, $u^*(D_j)=\sum_km'_{j,k}.D'_{j,k}$ are such that $m_j$ divides $m_{j,k}$ for any $j,k$. This cover is `orbifold-\'etale' if $m_j=m_{j,k}$, for every $j,k$.

In general, orbifold-ramified covers do not exist. An example is nevertheless the following: $u:\Bbb P^n\to (\Bbb P^n\vert \Delta)$, where $\Delta=\sum_{j=0}^{j=n}(1-\frac{1}{m_j}). H_j$, the $H_j$ being the $n+1$ coordinate hyperplanes.

We have then the following result ([Ca 07, th\'eor\`eme 6.33]): assume that $u:X'\to (X\vert \Delta)$ is an orbifold-ramified cover. Let $C'\subset X'$ be a rational curve not contained in $u^{-1}(Supp(\Delta))$. Then $C:=u(C')\subset X$ is a (`divisible') $\Delta$-rational curve. Conversely, if $C\subset X$ is a (`divisible') $\Delta$-rational curve, any component of $C':=u^{-1}(C)$ is rational if $u$ is orbifold-\'etale.
\end{example}

\begin{definition} The smooth pair $(X\vert \Delta)$ is uniruled (resp. rationally connected) iff any generic point of $X$ (resp. any generic pair of points of $X$) is contained in some $\Delta$-rational curve.
\end{definition}

\begin{example}\label{edrc}Let $(\Bbb P^n\vert \Delta)$, where $\Delta=\sum_{j=0}^{j=n}(1-\frac{1}{m_j}). H_j$, the $H_j$ being the $n+1$ coordinate hyperplanes be as in \ref{oec} above. Then $(\Bbb P^n\vert \Delta)$ is rationally connected, since $\Bbb P^n$ is rationally connected (in the usual sense).
\end{example}

We finish by a last conjecture, which extends to the orbifold context the uniruledness conjecture, which corresponds to the case $\Delta=0$:

\begin{conjecture}\label{corc} Let $(X\vert \Delta)$ be a smooth pair, then:

1. $\kappa(X\vert \Delta)=-\infty$ iff $(X\vert \Delta)$ is uniruled.

2. $\kappa_{++}(X\vert \Delta)=-\infty$ iff $(X\vert \Delta)$ is rationally connected.
\end{conjecture}

Of course, one would like to extend to the orbifold setting the  many facts  known when $\Delta=0$. But very few is known in this direction. For example, it is even unknown whether Fano smooth pairs are uniruled, this even in dimension $2$ (but the logarithmic case is then true, by [Keel-McKernan]). The case of Fano orbifolds is the decisive case for the solution of conjecture \ref{corc}, by the structure theorem \ref{tr^*}, which birationally expresses orbifolds with $\kappa_+=-\infty$ as towers of fibrations with Fano orbifold fibres. 

\begin{remark} We shall give a counting argument supporting the uniruledness of Fano smooth orbifolds, by showing that covering families of $\Delta$-nice rational curves (see \ref{exorc} for this notion) should exist in this situation.  Let indeed $g_0:\Bbb P^1\to X$ be a nonconstant map with $C:=g_0(\Bbb P^1)\subsetneq Supp(\Delta)$, going through a general point $a\in X$. Assume that $D_j.C=k_j.m_j$, for each $j\in J$, with $k_j$ an integer. The variety $Hom_a(\Bbb P^1,X)$ of such maps has at $g_0$ dimension $dim_{g_0}Hom_a(\Bbb P^1,X)\geq -K_X.C+3$. The number of conditions for $C$ to have order of contact at least $m_j$ at an (undetermined) point of $D_j$ lying on $C$ is equal to $m_j-1$. The total number of conditions for $g_0$ to be $\Delta$-nice' is thus $\sum_j k_j.(m_j-1)=\sum_j(1-\frac{1}{m_j}).k_j.m_j=\sum_j (1-\frac{1}{m_j}).D_j.C=\Delta.C$.

The expected dimension of the variety of such `$\Delta$-nice' rational curves through $a$ is thus at least $-(K_X+\Delta).C+3$, which thus remains positive after forgetting the $3$-dimensional space of parametrisations of $\Bbb P^1$, precisely when $(X\vert \Delta)$ is Fano.

\end{remark}

\begin{example}\label{etm} Let us consider the case when $X= \Bbb P^n, n\geq 2$, and when the support of $\Delta$ consists of $k$ hyperplanes $H_j$ in general position, with finite {\it integral} multiplicities $(m_0,\dots, m_{k-1})$, with $2\leq m_0\leq m_1\leq \dots\leq m_{k-1}$, in which case we shall just say that $(\Bbb P^n\vert \Delta)$ is of type $(m_0,\dots, m_{k-1})$. The condition that $(\Bbb P^n\vert \Delta)$ be Fano (resp. has trivial canonical bundle) is then just that $\sum_j(1-\frac{1}{m_j})<(n+1)$  (resp. that: $\sum_j(1-\frac{1}{m_j})=(n+1))$. The Fano condition is thus always satisfied when $k\leq (n+1)$. When $k\leq (n+1)$ it is not difficult to see (see \cite{Ca07}, and \ref{oec} above), that, for any finite set of points of $\Bbb P^n$, none of them lying on $\Delta$, there is an irreducible $\Delta$-rational curve\footnote{And even a `divisible' one, see \cite{Ca07} for this notion.} containing these points. See example \ref{k=n+1} below for a direct proof of their uniruledness. By contrast, when $k=(n+2)$, it is not known whether these Fano orbifolds are rationally connected. We shall give examples in which it can be shown by specific methods that they are, at least $\Delta$-uniruled (and covered by $\Delta$-elliptic curves when their canonical bundle is trivial). Observe that, for any $n\geq 2$, there is anyway only a finite number of $(n+2)$-tuples of integers $(m_0,\dots, m_{n+1})$such that $\sum_{j=0}^{n+1}(1-\frac{1}{m_j})\leq(n+1)$ (of course, provided that $\sum_{j=0}^{n}(1-\frac{1}{m_j})\leq n,$ see \ref{ss} below for this finiteness statement and examples).
\end{example}

\begin{example}\label{k=n+1} Assume first that $(\Bbb P^n\vert \Delta)$ is of type $(m_0,\dots, m_k)$, with $k\leq n$. Then the $H_j's$, for $j=1,\dots, k-1$ intersect in a projective space $P$ of dimension $n-(k-1)\geq 0$. Any projective line meeting $P$, but not contained in $P$, meets the support of $\Delta$ in at most $2$ points, and has finite $\Delta_g$ multiplicities $(m_0$ and $m_{k-1})$ there, and is thus $\Delta$-rational, and $(\Bbb P^n\vert \Delta)$ is thus uniruled by the family of lines through $P$.
\end{example}

\begin{example}\label{jneqn+1} Assume next that $(\Bbb P^n\vert \Delta)$ is of type $(m_0,\dots, m_{n+1})$, and that any smooth orbifold $(\Bbb P^{n-1}\vert \Delta')$ of type $(m_0,\dots,m_{n-2},m_{n-1}, m_{n+1})$ is $\Delta'$-uniruled (just $m_{n}$ has been omitted) if $\sum_{j\neq n}\frac{1}{m_j}>1$ . Then $(\Bbb P^n\vert \Delta)$ is $\Delta$-uniruled if $\sum_{j\neq n}\frac{1}{m_j}>1$. Indeed, the generic member of the pencil of hyperplanes $H_s$ containing $P:=H_{n}\cap H_{n+1}$ is naturally equipped with the orbifold structure $\Delta_s:=\sum_{j\neq n}(1-\frac{1}{m_j}).H_{s,j}$, where $H_{s,j}:=H_J\cap H_s$, so that $H_{s,n+1}=P$. Moreover, for each curve $g:C\to H_s$ birational onto its image, not contained in the union of the $H_{s,j}$, the orbifold divisor on $C$ computed from $(H_s\vert \Delta_s)$ and $(\Bbb P^n\vert \Delta)$ coincide (this is an immediate check). Assuming that $\sum_{j\neq n}\frac{1}{m_j}>1$, we deduce from the assumption made, that $H_s$ is $\Delta_s$-uniruled. Thus $(\Bbb P^n\vert \Delta)$ is uniruled, too.

(The same statement should hold, assuming that $\sum_{j\neq n+1}\frac{1}{m_j}>1$, but no obvious geometric construction seems to give this).
 \end{example}

\begin{example}\label{jneqn+1}  The first case when $k=n+2$ is when $X=\Bbb P^2$, and $\Delta=\sum_{j=0}^{j=3} (1-\frac{1}{m_j}). D_j$ is supported on $4$ lines in general position, of multiplicities $(m_0,m_1,m_2,m_3)$. Then $(X\vert \Delta)$ is Fano if and only if $\sum_j \frac{1}{m_j}>1$. This is the case when $(m_0,m_1,m_2,m_3)=(2,3,7,41)$, for example. It is then easy to show that a line is a $\Delta$-rational curve if and if only if it goes through $2$ of the $6$ double points of the union of the $4$ lines, so there are only $15$ such lines. But there is a one-dimensional family of conics which are $\Delta$-rational: the conics $C$ which are tangent to each of the $4$ lines. Indeed, for such a generic smooth conic $g:C\to \Bbb P^2$, $g$ being the incusion, the divisor $\Delta_g$ is supported on the $4$ points of tangency with multiplicities $(\frac{m_0}{2}, \frac{m_1}{2}, \frac{m_2}{2}, \frac{m_3}{2})$,  which is an orbifold rational curve, since: $-2+(1-\frac{2}{m_0})+(1-\frac{2}{m_1})+(1-\frac{2}{m_2})+(1-\frac{2}{m_3})=2.[1-(\frac{1}{m_0}+\frac{1}{m_1}+\frac{1}{m_2}+\frac{1}{m_3})]<0$. This Fano orbifold is thus indeed at least uniruled. By the very same argument, it is covered by $\Delta$-elliptic conics when its canonical bundle is trivial.\end{example}

\

We shall now partially extend the preceding example to higher dimensions.

\

\begin{example} Let us consider the case when $X= \Bbb P^n, n\geq 2$, and when the support of $\Delta$ consists of $(n+2)$ hyperplanes $H_j, j=0,1,\dots,(n+1)$ in general position (i.e: such that the intersection of any  $(n+1)$ of them is empty), with finite multiplicities $(m_0,\dots, m_{n+1})$, such that $\sum_j(1-\frac{1}{m_j})<(n+1)$, or equivalently, that $\sum_j(\frac{1}{m_j})>1$. We can, and shall, assume that the first $(n+1)$ hyperplanes $H_j, j=0,\dots,n$ are the coordinate hyperplanes of equations $X_j=0$, and that the last hyperplane $H_{n+1}$ has equation $X_0+\dots+X_n=0$, since all $(n+2)$-tuples of hyperplanes in general position are equivalent under homographies.

{\bf Remark:} Assume that $C\cong \Bbb P^1$ is the normalisation of an irreducible rational curve of degree $d$ on $\Bbb P^n$, meeting each of the hyperplanes $H_j$ in a single point $a_j$, thus with contact order $d$. Let us try to determine a condition on $d$ which implies that $C$ is $\Delta$-rational. The orbifold multiplicity at such a point $a_j\in C$ is thus: $\frac{m_j}{d}$. The corresponding orbifold divisor on $C$ thus consists of $(n+2)$ points with multiplicities $\frac{m_j}{d}$, and this curve is thus $\Delta$-rational if and only if: $\sum_j\frac{d}{m_j}>(n+2)-2=n$, or equivalently, if: $\sum_j\frac{1}{m_j}>\frac{n}{d}$. Since, by assumption, $\sum_j \frac{1}{m_j}>1$, this condition is realised as soon as $\frac{n}{d}\leq 1$. We thus may choose: $d=n$, and $C$ to be a rational normal curve of degree $n$. Notice, however, that we did not take into account the fact that the multiplicity at $a_j$ is taken to be $1$, and not $\frac{m_j}{d}<1$ if $m_j<d$. This is discussed in remark \ref{rvrc} below.

Recall that a normal rational curve of degree $n$ on $\Bbb P^n$ is parametrically given by a map: $P(t)=(P_0(t):\dots:P_n(t))$, where the $P_j(t)' s$ are linearly independent polynomials of degree $n$. All such curves are equivalent under the natural action of $\Bbb PGl(n+1,\Bbb C)$ $\square$

\begin{theorem}\label{tnrc} For any set of $(n+2)$ hyperplanes $H_j$ in general position on $\Bbb P^n$, and for any $p=(p_0:\dots:p_n)\in \Bbb P^n$ generic, there exists a rational normal curve $C$ of degree n on $\Bbb P^n$, which goes through $p$, and meets each of the $(n+2)$ hyperplanes $H_j$ in exactly one point, which lies on the smooth part of the union of the $H_j's$.

Remark however that such a curve is `virtually $\Delta$-rational', but not always $\Delta$-rational if $n\geq 3$ (see remark \ref{rvrc} below). \end{theorem}
\end{example}

{\bf Proof:} The condition that $C$ meets each of the coordinate hyperplanes $H_j$ in a single point thus reads as $C$ being given by a parametric representation of the form: $P(t)=(b_0.(t+a_0)^n,\dots,b_n.(t+a_n)^n)$ for nonzero and pairwise distinct complex numbers $a_j$, and nonzero complex numbers $b_j$, $j=0,\dots,n$.

The condition that $C$ meets the $(n+2)$-th hyperplane $H_{n+1}$ translates as the equation: $b_0.(t+a_0)^n+\dots+b_n.(t+a_n)^n=b.(t+a)^n$ for nonzero numbers $b=\sum b_j,a$, with $a$ being distinct from all of the preceding $a_j's$.

Because we also want $C$ to go through the point $p=(p_0,\dots,p_n)$, at time $t=\infty$ say, we get the additional conditions: $b_j=p_j$.

In our situation, the $p_j's$ are given, and the numbers $a,b_j,a_j's$ are to be determined from these algebraic equations. These equations are however of high degree and difficult to solve. We shall proceed differently: solve for the $b_j's$ assuming the $a, a_j's $ to be given, because this is simply a linear system. And then show that the $a, a_j's$ exist, for a generic choice of the $p_j's$, after projectivisation. This is sufficient to imply the result, by chosing the $p_j's$ generically.

Writing the $(n+1)$ coefficients of the two polynomials in $t$ on the two side of the equation $b_0.(t+a_0)^n+\dots+b_n.(t+a_n)^n=b.(t+a)^n$ written above permits to rewrite this equation as a linear system in the $p_j's$. In matrix form, it reads as: $A.v=b.w$, with $A$ the square complex matrix of size $(n+1)$ given below, $^t v:=(p_0,\dots,p_n)\in (\Bbb C^*)^{n+1}$, and $^tw:=(1,a,a^2,\dots,a^n)\in (\Bbb C^*)^n$.

\

The matrix $A$ is:

\

$A:=\left(
\begin{array}{ccccc}
1&1&\dots& 1&1\\
a_0 &a_1&a_2&\dots&a_n \\
a_0^2&a_1^2&a_2^2&\dots&a_n^2\\

\dots&\dots &\dots&\dots\\
a_0^n&a_1^n&a_2^n&\dots&a_n^n\\
\end{array} \right)$

\

Now Cram\' er's rule permits to express the $p_j=v_j$ as a quotient of two determinants of Vandermonde type, in terms of the $b,a,a_j's$ supposed to be given. We finally get, after some simplification: $v_j:=p_j=b. \Pi_{h\neq j}\big(\frac{a-a_h}{a_j-a_h}\big)$, $h$ running from $0$ to $n$, avoiding $j$. Observe that the right-hand side is invariant by translation (i.e: takes the same value when we replace $a,a_j$ by $a+t,a_j+t$, for any $t\in \Bbb C)$.

We are thus reduced, by the above translation invariance, which permits to take $a=0$, to prove that the following rational map $\Phi$ is dominant.

The map $\Phi:\Bbb C^{n+1}\dasharrow \Bbb P^{n}$  is defined by $\Phi(a_0,\dots,a_n):=[x_0:\dots:x_n]$, with $x_j:=\Pi_{h\neq j}\big(\frac{a_h}{a_h-a_j}\big)=\Pi_{h\neq j}\big((1-(\frac{y_h}{y_j}))^{-1}\big)=\Psi(y_1,\dots,y_n)$, replacing the $(n+1)$ variables $a_0,\dots,a_n$ by the $n$ variables $y_j:=\frac{a_0}{a_j}$, for $j=1,\dots, n$. We denote also for notational simplication: $y_0:=\frac{a_0}{a_0}=1$.

Thus $\Phi$ is dominant if the determinant $det (J)$ of the Jacobian matrix of the logarithmic derivatives of the $n$ functions $u_j:=\frac{x_j}{x_0}$ does not vanish at some point where the functions $u_j=u_j(y_1,\dots,y_n)$ are regular and non zero.

Let us first rewrite, after some simplification: $u_j:=-y_j^n.\Pi_{h>0, h\neq j}.\Big(\frac{1-y_h}{y_j-y_h}\Big)$

\

Taking logarithmic derivatives now shows that:

\

$\bullet$ $\frac{y_k}{u_j}.\frac{\partial u_j}{\partial y_k}=-\frac{(1-y_j)}{(1-y_k)(1-\frac{y_j}{y_k})})$ if: $1\leq j\neq k\geq 1$, and that: 

$\bullet$ $\frac{y_j}{u_j}.\frac{\partial u_j}{\partial y_j}=n-\sum_{h\neq j}\Big((1-\frac{y_h}{y_j})^{-1}\Big).$

\

Let us now choose $M>0$, and the $y_j's$ all nonzero in such a way that: $\vert y_j\vert>M.\vert y_{j-1}\vert$ for $j=2,\dots, n$, with $M.\vert y_n\vert<1$. As $M$ goes to $+\infty$, one easily checks, using the equalities above, and since $\frac{y_j}{y_k}$ goes to $0$ if $j<k$ and to $\infty$ if $j>k$, that:

1. $\frac{y_k}{u_j}.\frac{\partial u_j}{\partial y_k}$ goes to $-1$ if $j< k$, and goes to $0$ if $j>k$, while:

2. $\frac{y_j}{u_j}.\frac{\partial u_j}{\partial y_j}$ goes to $n-(j-1)$, for any $j=1,\dots,n$.

Thus $det (J). y_1.\dots.y_n$ goes to $det (J_0)$, where $J_0$ is the matrix with coefficients $-1$ below the diagonal, with coefficients $0$ above the diagonal, and coefficient $n-(j-1)$ on the diagonal, at the intersection of the $j$-th line and $j$-th row, for $j=1,\dots, n$. Since $det (J_0)=n!\neq 0$, $det (J)\neq 0$ when the $y_j's$ satisfy the above inequalities for $M$ sufficiently large , which implies the desired assertion that $\Phi$ is dominant $\square$

\

From the theorem \ref{tnrc} above we shall now deduce that smooth orbifolds are covered by `virtual' $\Delta$-rational or elliptic curves, according to whether they are Fano, or have trivial canonical bundle. We define first these `virtual' notions.

\

\begin{definition}\label{dvrc} Let $(X\vert \Delta)$ be a smooth orbifold, with $\Delta:=\sum_{j\in J}(1-\frac{1}{m_j}).D_j$. Let $C$ be a smooth connected projective curve. A map $g:C\to (X\vert \Delta)$ is a `virtual' $\Delta$-rational (resp. a `virtual' $\Delta$-elliptic) curve if:

1. It is birational onto its image, which is not contained in $Supp(\Delta)$.

2. $deg(K_{C}+\Delta^*_g)<0$ (resp. $deg(K_{C}+\Delta^*_g)=0)$, where $\Delta^*_g$ is the orbifold divisor on $\Bbb P^1$ which assigns to any $a\in \Bbb P^1$ the multiplicity $1$ if $g(a)\notin Supp(\Delta)$, and otherwise the multiplicity $m_g(a):=max_{j\in J(a)}\{\frac{m_j}{t_{j,a}}\}$, where $J(a), t_{j,a}$ are defined as in \ref{dorc}.
\end{definition}

\begin{example}\label{rvrc} Assume now that $(X\vert \Delta)$ is of type $(m_0,\dots, m_{n+1})$ on $X=\Bbb P^n$, in the sense of example \ref{etm}, and that $C$ is a rational normal curve of degree $n$ meeting each of the hyperplanes $H_j$ in one point with contact order $d$, as in theorem \ref{tnrc}. Then $C$ is virtually $\Delta$-rational or elliptic, according to whether $(X\vert \Delta)$ is Fano, or has trivial canonical bundle (see corollary \ref{crnc} below). In general, $C$ will not be $\Delta$-rational, or elliptic, unless $n=2$. See proposition \ref{vvsnv} below.
\end{example}

\begin{corollary}\label{crnc} Let $(\Bbb P^n\vert \Delta)$ be a smooth orbifold, the support of $\Delta$ consisting of the union of $(n+2)$ hyperplanes $H_j$ in general position. The generic point of $\Bbb P^n$ is then contained in a normal rational curve $C$ of degree $n$ which meets each $H_j$ in a single point, and such a curve is virtually $\Delta$-rational (resp. $\Delta$-elliptic) if $(\Bbb P^n\vert \Delta)$ is Fano (resp. has trivial canonical bundle).

Moreover, if either $n=2$, or $m_1\geq n$, or more generally, if: $\sum_{j=0}^{j=n+1}\frac{1}{m_j^*}>1$ (resp. if: $\sum_{j=0}^{j=n+0}\frac{1}{m_j^*}=1)$, then $C$ is also $\Delta$-rational (resp. $\Delta$-elliptic), where $m_j^*:=max\{m_j,n\}$.
\end{corollary}

{\bf Proof:} The first assertion follows in the Fano case from theorem \ref{tnrc} and the remark made before its statement. The same computation works in the trivial canonical bundle case: let $m_j$ be the multiplicities of the $H_j's$. The orbifold
$(\Bbb P^n\vert \Delta)$ has trivial canonical bundle if and only if:  $\sum_j\frac{1}{m_j}=1.$ The orbifold multiplicity at a point $a_j\in C$ is thus: $\frac{m_j}{n}$. The corresponding orbifold divisor on $C$ thus consists of $(n+2)$ points with multiplicities $\frac{m_j}{n}$, and this curve is thus $\Delta$-elliptic if and only if: $\sum_j\frac{n}{m_j}=(n+2)-2=n$, which holds true, since, by assumption, $\sum_j \frac{1}{m_j}=1$

Let us now check the second assertion. If $m_1\geq n$, then $m_j\geq n, \forall j$, so that $\frac{m_j}{n}\geq 1,\forall j$, and so no $max$ is needed to compute $\Delta_g$ in definition \ref{dorc}. The conclusion thus follows. Notice that $m_1\geq 2$, so that the conclusion always holds when $n=2$. 

Now if , for some $1\leq j\leq n+2$, $m_j<n$, taking $max\{1,\frac{m_j}{t_{j,a}}=\frac{m_j}{n}\}$ amounts to replacing $m_j$ by $m_j^*$, and so $\Delta^*_g$ is simply $\Delta_g$ computed for $(X\vert \Delta^*)$ instead of $(X\vert \Delta)$, with $\Delta^*:=\sum_{j=0}^{j=n+1}(1-\frac{1}{m_j^*}).H_j$, which implies the conclusion by the first part $\square$

We thus see that the consideration of rational normal curves of degree $n$ permits to show the uniruledness of some Fano smooth orbifolds of type $(m_0,\dots, m_{n+1})$, and of all if $n=2$. Unfortunately, when $n\geq 3$, these are, by far, not all Fano orbifolds of this type. We shall make now more precise which are the Fano smooth orbifolds of dimension $3=n$ which can be shown to be uniruled by this method, and which are not. First collecting the results of \ref{tnrc}, \ref{jneqn+1}, and \ref{crnc}, we get the assertions 1,2,3 below:

\begin{proposition}\label{vvsnv} Let $(\Bbb P^n\vert \Delta)$ be a smooth Fano orbifold of type $(m_0,\dots, m_{n+1})$. Then $(\Bbb P^n\vert \Delta)$ is uniruled, unless (maybe) if the following $3$ conditions are realised:

1. $\sum_j\frac{1}{m_j}>1$

2. $\sum_{j\neq (n)}\frac{1}{m_j}\leq1$

3. $\sum_j\frac{1}{m^*_j}\leq1$, where $m_j^*:=max\{m_j,n\}$, if $n=3$.

4. If $n=3$, then $(\Bbb P^3\vert \Delta)$ is uniruled, unless possibly if $m_0=2, m_1\geq 3$, and $m_2\geq 4$.

\end{proposition}

{\bf Proof:} We have only to prove the assertion $4$. If $m_0\geq 3$, then $m_j^*=m_j,\forall j$, and so \ref{vvsnv}.3 above contradicts  \ref{vvsnv}.1. We thus assume that $m_0=2$ in the sequel. In the same way, if $m_1=2$, or if $m_1=m_2=3$, then the sum of the first $2$ or $3$ terms of $\sum_j\frac{1}{m_j}$ is at least $1$, contradicting \ref{vvsnv}.2 above. Thus $m_2\geq 4$.
$\square$

When $n=3$, the `Fano' types $(m_0,m_1,m_2,m_3,m_4)$ for which the uniruledness can (or cannot be) proved by the preceding method can be (lengthily) listed. We shall give in \ref{tec} two extremal cases when $n=3$, for which the uniruledness $(\Bbb P^3\vert \Delta)$ cannot be proved by the preceding method.

\begin{definition}\label{dgss} Let $(2\leq a_1\leq a_2\leq \dots\leq a_k)$ be a finite sequence of positive integers such that $\sum_{j=1}^{j=k}\frac{1}{a_j}=1-\frac{1}{b}$, for some integer $b\geq 2$. 

Define inductively for $s\geq 1$: $a_{k+1}:=b+1, a_{k+s+1}:=a_{k+s}.(a_{k+s}+1)+1$.

It is then immediate to show inductively that: $\sum_{j=1}^{j=k+s}\frac{1}{a_j}=1-\frac{1}{a_{k+s+1}-1}$, for any $s\geq 1$.

Thus $\sum_{j=1}^{j=k+s}\frac{1}{a_j}+\frac{1}{a_{k+s+1}-2}=1+\frac{1}{(a_{k+s+1}-2).(a_{k+s+1}-1)}>1$, for any $s\geq 0$. This will give us examples of Fano types on projective spaces.
\end{definition}

\

\begin{example}\label{egss} We look at special cases of the preceding sequences. They provide examples of types of Fano orbifolds which are uniruled, by \ref{tnrc} and its corollaries.

\

Let $k\geq 1$, and choose $a_1=\dots=a_k=k+1$, so that  $b=k+1$. We then get: $a_{k+1}=k+2,a_{k+2}=(k+2)(k+3)+1, a_{k+3}=a_{k+2}.(a_{k+2}+1)+1,\dots$.

When $k=1$, we get: $d_1=2,d_2=3,d_3=7,d_4=6.7+1=43, d_4=42.43+1=1807,\dots$  (denoting $a_j=d_j$, in this case).

When $k=2$, we get: $t_1=t_2=3,t_3=4,t_4=13, t_4=157,\dots$ (denoting $a_j=t_j$, in this case).

When $k=4$, the sequence is: $q_1=q_2=q_3=4,q_4=5,q_5=21,q_6=421,\dots$  (denoting $a_j=q_j$, in this case).

\

When $k=n-1$, the types $(m_1=n,\dots,m_{n-1}=n, m_n=(n+1), m_{n+1}=n.(n+1)+1,m_{n+2}=n.(n+1).(n^2+n+1)-2$ are the types of  Fano orbifolds $(\Bbb P^n\vert \Delta)$ with orbifold divisor supported on the union of $(n+2)$ hyperplanes in general position, and all such orbifolds are uniruled, by theorem \ref{tnrc} and its corollaries.

\

Specific examples are thus: $n=3$, and type $(3,3,4, 13, 155)$, or: $n=4$ and type $(4,4,4,5,21,419)$.
\end{example}

\

\begin{proposition}\label{ss} Let $N\geq 1$ be an integer. There exists a bound $B_N<1$ such that if $2\leq a_1\leq \dots \leq a_N$ is a finite sequence of integers, and if $A:=\sum_{j=1}^{j=N}a_j<1$, then : $A\leq B_N$.
\end{proposition}

{\bf Proof:} We assume the existence of $B_N<1$ and are going to establish inductively  the existence of $B_{N+1}<1$. It is plain that $B_1=\frac{1}{2}$ exists. Assume now that $A+\frac{1}{a_{N+1}}:=A'=\sum_{j=1}^{j=N+1}\frac{1}{a_j}<1$. We can increase strictly $A'$, preserving this last inequality, by replacing $a_{N+1}$ by $a_{N+1}-1$ (and possibly reordering the terms, in case $a_N=a_{N+1})$, unless:
$B_N+\frac{1}{a_{N+1}-1}\geq A+\frac{1}{a_{N+1}-1}\geq 1$. 
 Thus, if $A'$ cannot be increased in this way, as we may assume, we have: $\frac{1}{a_{N+1}-1}\geq (1-B_N)$, and $a_{N+1}\leq 1+\frac{1}{1-B_N}$. There are thus only finitely many values for all $a_j's$, and there exists some $B_{N+1}<1$ such that $A'\leq B_{N+1}$ $\square$

\begin{example}\label{tec} The types $(m_0,m_1,m_2,m_3,m_4)$ for which $(\Bbb P^3\vert \Delta)$ is Fano, but which cannot be proved to be uniruled by the preceding method satisfy in particular: $m_0=2, 3\leq m_1\leq 7, 4\leq m_2$ (this is easy, using \ref{vvsnv}). There are two main cases:

a. $\sum_{0\leq j\leq 3}\frac{1}{m_j}<1$. There are only finitely many of them, and then $m_3\leq d_4=43$, as may be shown using \ref{ss}. Then also $m_4< d_5-2$. 

A typical example is $(2,3,7,43, 1805)$.

b. $\sum_{0\leq j\leq 3}\frac{1}{m_j}\geq1$. There are only finitely many such $4$-uplets, since $\sum_{0\leq j\leq 2}\frac{1}{m_j}<1$, as follows from \ref{vvsnv}.3, and so $m_3\leq 42=d_4-1$. However, for any sufficiently large $m_4$, the given `type'  satisfies the inequalities of \ref{vvsnv}. 

Typical examples are $(2,3,7,42,m_4)$, with any $m_4\geq 42$. 
\end{example}

\begin{remark} To show the uniruledness of Fano orbifolds of type $(2, 3,m_2, m_3,m_4)$ on $\Bbb P^3$, with $m_2\geq 6$, it were sufficient to show, for any $5$-uple of hyperplanes $H_j,j=0,\dots, 4$ of $\Bbb P^3$, the existence of a rational curve $C$ of degree $6$ tangent in $3$ points to $H_0$, meeting $H_1$ in two distinct points with order of contact of order $3$, and meeting each one of the three remaining $H_j's$ in one single point with order of contact $6$. Indeed: the resulting $C$ multiplicities were then: $1$ for the first $5$ points, and  $\frac{m_j}{6}$, for $j=3,4,5$, and the last $3$ points. This were a $\Delta$-rational curve, since: $\sum_{2\leq j\leq 4}\frac{6}{m_j}= 6.[(\sum_{0\leq j\leq 4}\frac{6}{m_j})-(\frac{1}{2}+\frac{1}{3})]>6.(1-\frac{5}{6})=1$.

This construction appears to be similar to the one made above for the rational normal curves of degree $3$. A simple dimension count shows that that such curves should exist. The general case $n\geq 4$ seems to require other ideas and techniques, however
\end{remark}

\

\section{\bf Pseudoeffectivity of the relative canonical bundle}

\

In this section, we shall prove by a relative Bend-and-Break technique, that the relative canonical bundle of a fibration is pseudoeffective if the generic fibre is not uniruled, and derive from it a weak version of Viehweg's weak positivity for the direct images of pluricanonical sheaves. Although the results are weaker than known ones, the method of proof is so straightforward, and possibly susceptible of further developments, that it seemed worth being written. Combined with Hodge-theoretic arguments, it  might indeed permit to easily obtain stronger versions, closer to Viehweg's results.

\begin{theorem}\label{tprk} Let $f:X\to Y$ be a fibration, with $X,Y$ smooth projective connected. Assume that some smooth fibre $X_y$ of $f$ is not uniruled. Then $K_{X/Y}$ is pseudo-effective.
\end{theorem}

{\bf Proof:} Assume that $K_{X/Y}$ is not pseudo-effective. By [BDPP 04], there exists an algebraic $X$-covering family $(C_t)_{t\in T}$ of curves such that $-K_{X/Y}.C_t>0$. We can assume that the generic curve of this family is not rational, since it may be obtained as the direct image of a complete intersection of very ample divisors on some blow-up of $X$. Let $a\in X$ be a general point, lying on some smooth fibre $X_y$ of $f$, and also on some nonrational irreducible member $C$ of the family $(C_t)_{t\in T}$. Since $-K_{X/Y}.C>0$, there exists, by [De 01], Proposition 3.11, p. 70,  and the dimension estimate (2.4), p. 47, using the now standard Mori reductions to characteristic $p>0$, a rational curve $R\subset X_y$ passing through $a$. Since $X_y$ is not uniruled, choosing $a$ not lying on any rational curve contained in $X_y$, we get a contradiction to our initial assumption that $K_{X/Y}$ is not pseudoeffective $\square$

\

Remark that the relative Bend-and-Break lemma used above is the same as the one used in [Ca 92] and [KMM 92] to show the rational chain-connectedness of Fano manifolds. 

\

We now deduce from the preceding theorem \ref{tprk} a proof of proposition \ref{pabiad} when $\Delta=0$.

\begin{corollary}\label{cprk} Let $f:X\to Y$ be a fibration, with $X,Y$ smooth, $X$ projective, and let $X_y$ be a generic fibre of $f$. Then:

1. Assume $\nu(X,K_X)=\kappa(X)$. Then $\kappa(X)\geq \kappa(X_y)+\kappa(Y)$.


2. If $X_y$ and $Y$ are of general type, so is $X$.
\end{corollary}

{\bf Proof:} 1. We can assume that $\kappa(Y)\geq 0$. Let $A$ be any ample divisor on $X$. Then $m.K_Y$, and also $m.K_{X/Y}+A$ are effective, for some $m>>0$, by theorem \ref{tprk} above. Thus $mK_{X}+A$, and so $N.K_X$ is effective, too, for some $N>>0$, since $\nu(X,K_X)=\kappa(X)$, by our assumption. The claim then follows from the arguments of \cite{AC}, lemma 2.4, p. 516, for example.

2. This follows from the lemma \ref{lb} below, applied to $P:=K_{X/Y}$ and $D:=K_Y$$\square$

\begin{lemma}\label{lb} Let $f:X\to Y$ be a fibration. Let $P$ be a pseudo-effective line bundle on $X$ which is $f$-big (i.e: big on the generic fibre $X_y$ of $X$). 

Then: for any big $\bQ$-divisor $D$ on $Y$,  $P+\epsilon.f^*(D)$ is also big.\end{lemma}

{\bf Proof:} $P+f^*(m.D)$ is big on $X$ for some $m>>0$, by the assumption of relative bigness. Since $P$ is pseudo-effective, $(N-1)P+(P+m.f^*(D))=N.(P+\frac{m}{N}.f^*(D))$ is also big, for any $N\geq 1$. Choosing $N\geq m$ gives the claim $\square$

\begin{remark} 1. The proof of theorem \ref{tprk}does not seem to be able to give the weak positivity statement given by Viehweg's theorem. It also does not apply (directly at least) to the `orbifold' context of pairs $(X\vert \Delta)$. In this respect, it is much weaker. 

2. There is one point for which it is, however, more flexible: it does not need the effectivity of $K_{X_y}$, and its proof also directly gives information on multiples $m.K_{X/Y}$, contrary to Viehweg's proof which requires two steps: dealing first with $m=1$, and then with arbitrary $m's$.

3. The corollary \ref{cprk}(2) is known in a much more stronger version, by \cite{K88}, which proves $C_{n,m}^+$ when $X_y$ is of general type.
\end{remark}

\section{\bf Families of canonically polarised manifolds.}

Roughly stated, a generalisation by Viehweg of a conjecture of Shafarevich states that `the moduli space of canonically polarised manifolds has components of log-general type'. The initial formulation was that if $f:X\to B$ is an algebraic smooth family of canonically polarised manifolds parametrised by a quasi-projective manifold $B$ having generically a `variation' (i.e: a Kodaira-Spencer map) of maximal rank, then $B$ is of log-general type, considering a smooth compactification $B=Y-D$, such that $(Y\vert D)$ is smooth, with $D$ reduced and of snc.

This has been shown in low dimension and various formulations by Viehweg-Zuo, Kebekus-Kovacs, Jabbusch-Kebekus (see \cite{VZ} and \cite{JK} for the appropriate references).

The natural formulation of this conjecture seems to be:

\begin{conjecture}\label{ci} (`isotriviality conjecture') Let $f:X\to B$ be as above, assume that $B$ is special (i.e: that so is any smooth compactification $(Y\vert D)$ as above. Then $f$ is isotrivial (i.e: all fibres of $f$ are isomorphic).
\end{conjecture}

This implies Viehweg's conjecture, since the moduli map (for arbitrary families $f:X\to B$) then factors through the core of $(Y\vert D)$, which is of log-general type.

\begin{theorem}\label{tlmmfcp} The Isotriviality Conjecture follows from conjectures $C_{n,m}^{orb}$ and \ref{c5}.

\end{theorem}

\begin{remark} The isotriviality conjecture is thus reduced to standard conjectures of birational geometry\footnote{It was stated in [Ca 09] that the isotriviality conjecture follows from conjecture \ref{ck++}. But it is only true that it follows from the arguments used to deduce conjecture \ref{ck++} from conjectures $C_{n,m}^{orb}$ and \ref{c5}. The proof given in the present text simplifies the arguments given in \cite{Ca 09}.}.\end{remark}

\

{\bf Sketch of proof:} The proof rests essentially on the construction by Viehweg-Zuo of a line bundle $L\subset Sym^m(\Omega^1_Y(log D))$ such that $\kappa(Y,L)=Var(f)$ (see [V-Z 02]). Because, as a consequence of $C_{n,m}^{orb}$ and theorem \ref{tc=jr}  we have a canonical decomposition $c=(J\circ r^*)^n$ which is the constant map (since $B=Y-D$ is special), it is sufficient to show the result when either $\kappa(B)=0$, or $\kappa_+(B)=-\infty$. In the first (resp. second) case, we have, by Conjecture \ref{c5} (1), the existence of a sequence of divisorial contractions and flips $s:(X\vert \Delta)\dasharrow (X'\vert \Delta')$ such that $c_1(X'\vert \Delta')=0$ (resp. such that a non-trivial fibration $f:(Y\vert D)\to Z$ with Fano orbifold fibres $(Y'_z\vert \Delta'_z)$ exists). In the first case, we directly conclude that $\kappa(Y,L)\leq 0$, which implies that $Var(f)=0$, as claimed. In the second case, we are reduced, by the equality $r^*=r^n$ of theorem \ref{tr^*} to the case where $(Y\vert D)$ is Fano. Considering, as in the proof of lemma \ref{fanok++} above, a new orbifold divisor $\Delta^+=D+\frac{1}{m.n}.H>D$, with $H\in \vert -m.n.(K_Y+D)\vert$, so that $c_1(Y\vert \Delta^+)=0$, we conclude $\kappa(Y\vert D,L)=-\infty$, since $\kappa(Y\vert \Delta^+)\leq 0$, so that the family is isotrivial on these fibres. $\square$

\

\end{document}